\documentclass[a4paper,12pt,legno]{article}

\usepackage{amsmath}
\usepackage{amsfonts}
\usepackage{amsthm}
\usepackage{mathrsfs}

\newtheorem{thm}{Theorem}

\newtheorem{lem}[]{Lemma}
\newtheorem{cor}{Corollary}

\newtheorem{rem}[]{Remark}

\newcommand{\Rset}{\mathbb{R}}
\newcommand{\Cset}{\mathbb{C}}

\newcommand{\La}{\mathcal{L}}

\newcommand{\al}{\alpha}
\newcommand{\be}{\beta}

\newcommand{\bl}{\bigl(}
\newcommand{\br}{\bigr)}

\newcommand{\wh}{^{\widehat{}}}

\begin{document}

\title{A random integral calculus on generalized s-selfdecomposable probability measures\footnote{Research funded by Narodowe Centrum Nauki (NCN)
Dec2011/01/B/ST1/01257}}

\author{ Zbigniew J. Jurek}

\date{June 22, 2013 }

\maketitle
\begin{quote} \textbf{Abstract.} The class $\mathcal{U}_{\beta}$ of
generalized s-selfdecomposable probability distributions can be
viewed as an image, via the random integral mapping
$\mathcal{J}^{\beta}$, of the class $ID$ of all infinitely divisible
measures. We prove that a composition of the mappings
$\mathcal{J}^{\beta_1}, \mathcal{J}^{\beta_2}, ...,
\mathcal{J}^{\beta_n}, \,\be_1>0, ..., \be_n>0,\,$ is again a random
integral but with a new deterministic inner time.
Moreover, some
elementary formulas concerning the distributions of products  of
powers of independent uniformly distributed random variables are
established.

\emph{Mathematics Subject Classifications}(2000): Primary 60F05 ,
60E07, 60B11; Secondary 60H05, 60B10.

\medskip
\emph{Key words and phrases:} Class $\mathcal{U}_{\be}$
distributions; generalized s-selfdecomposable distributions;
infinite divisibility; L\'evy-Khintchine formula; L\'evy process;
random integral; uniform distributions; Euclidean space; Banach
space.

\emph{Abbreviated title:} On generalized s-selfdecomposable measures

\end{quote}

\newpage
Let a real $\beta$ and infinitely divisible probability measures
$\nu_j$ (on Euclidean space $\Rset^d$ or a Banach space $E$) be such
that
\[
T_{1/n}(\nu_1\ast\nu_2 \ast...\ast \nu_n)^{\ast n^{- \beta}}
\ast\delta_{x_n}\Rightarrow   \mu \,,  \ \ \mbox{as} \ \
n\to\infty\,,\ \ \ (\star)
\]
for some deterministic shifts $x_n$, where $T_c$ denotes the
dilation (multiplication) by $c>0$. Then we write
$\mu\in\mathcal{U}_{\beta}$ and call $\mu$ \emph{a generalized
s-selfdecomposable measure}. In a series of papers Jurek (1988,
1989), Jurek and Schreiber (1992) it was proved, among others, that
for not degenerate $\mu$ in $(\star)$ we must have $\beta \ge -2$
and that $\mathcal{U}_{\beta}$  ($\beta \ge -2$) form an increasing
family of convolution semigroups that $"$almost$"$ exhaust the whole
class $ID$ of all infinitely divisible measures.

In a recent paper by James and Zhang (2011) generalized s-
selfdecomposable distributions from $\mathcal{U}_{\beta}$
($\beta>0$) were used for price models that exhibit volatility
clustering.

For the purpose of this paper the most crucial is the fact that
generalized s-selfdecomposable measures admit the \emph{random
integral representation} (1), (see below), which is a particular
case of the following  representation \emph{
\begin{multline*}
\qquad \qquad \qquad  \mu\equiv I^{h,r}_{(a,b]}(\nu):=\mathcal{L}\big( \int_I h(t)dY_{\nu}(r(t))\big),    \ \ \  \ \ \ \  (\star \star) \\
\mbox{where} \ I=(a,b] \subset \Rset^+, \
h:\Rset^+\to \Rset, \ \ Y_{\nu}(\cdot) \ \mbox{is a L\'evy process and}, \\
r:\Rset^+ \to \Rset^+  \mbox{is a monotone function (deterministic
time change in $Y_{\nu}$) .}
\end{multline*}}
In fact, it was shown that many  classes of limit laws can be
described as collections of probability distributions of random
integrals of the form $(\star \star)$ for suitable chosen parameters
$h,\,r$ and $I$ (possibly a half-line). Later on, this led to the
conjecture that \emph{all classes of limit laws} should admit random
integral representation; cf. a survey article Jurek (2011) and see
the Conjecture on www.math.uni.wroc.pl/$\sim$zjjurek \footnote{It
might be of an interest to recall here that S. D. Chatterji's
\underline{subsequence principle} claiming that: \ \emph{Given a
limit theorem for independent identically distributed random
variables under certain moment conditions, there exists an analogous
theorem such that an arbitrary-dependent sequence (under the same
moment conditions) always contains a subsequence satisfying this
analogous theorem} was proved by David J. Aldous (1977). Although,
we do not expect that the above Conjecture and Chatterji's
subsequence principle are mathematically related, however, one may
see a $"$philosophical$"$ relation between those two.}.

One may hope that this note will lead to establishing of \emph{a
calculus on random integral mappings $I^{h,r}_{(a,b]}$}, and their
domains of definition $\mathcal{D}^{h,\,r}_{(a,b]}$, analogous to
that of the linear operator calculus in functional analysis.

In this paper we will prove that the class of the integral mappings
(1), for $\beta>0$, is closed under compositions, that is, their
compositions are of the form $(\star \star)$ with the properly
chosen time change $r$; (Theorem 1). As an auxiliary result we found
a decomposition of number 1 as a sum of products of complex
fractions; (Lemma 1). Also compositions of the mappings (1) are
described in terms of the L\'evy-Khintchine triples; (Theorem 2).
Auxiliary Lemma 3 give probability distribution functions (p.d.f.)
of products of powers of independent uniformly distributed random
variables as linear combinations of other p.d.f.

\medskip
 \textbf{1. Introduction and main results.}

\emph{The results here are given for random vectors in $\Rset^d$ .
However, proofs are such that they are valid for infinite
dimensional separable Banach spaces $E$ when ones replaces a scalar
product by the bilinear form on the product space $E^{'}\times E$,
where $E^{'}$ denotes the dual space; see Araujo and Gine (1980) or
Linde (1983) or Ledoux and Talagrand (1991). Of course,
$(\Rset^d)^{'}=\Rset^d$.}

\medskip
Throughout the paper $\mathcal{L}(X)$ will denote the probability
distribution of an $\Rset^d$-valued random vector $X$; (or a real
separable Banach space E-valued random element $X$ if the Reader is
interested in that generality). Similarly, by $Y_{\nu}(t), t\ge0,$
we will denote an $\Rset^d$- valued (or an E-valued) L\'evy process
such that $\mathcal{L}(Y_{\nu}(1))=\nu$. Recall that by a L\'evy
process we mean a process with stationary independent increments,
starting from zero, and with paths that are continuous from the
right and with finite left-hand limits (that is, cadlag paths). Of
course, $\nu\in ID$, where $ID$ stands for all
\emph{\underline{i}nfinitely \underline{d}ivisible} measures on
$\Rset^d$ (or on a Banach space $E$).

For $\beta>0$ and a L\'evy process $Y_{\nu}(t), t\ge0$, we define
mappings
\begin{equation}
\mathcal{J}^{\{\beta\}}(\nu)\equiv
\mathcal{J}^{\beta}(\nu):\,=\La\bl\int_{(0,1]}
t^{1/\beta}\;dY_{\nu}(t)\br=\La\bl\int_{(0,1]}
t\;dY_{\nu}(t^{\beta})\br
\end{equation}
and  the classes $\mathcal{U}_{\beta}:\,=\mathcal{J}^{\beta}(ID)$.
To the distributions from $\mathcal{U}_{\beta}$ we refer to as
\emph{generalized s-selfdecomposable distributions}.

These classes of probability measures were originally defined as
limiting distributions in some schemes of summations; cf. Jurek
(1988 and 1989). In particular, the class
$\mathcal{U}\equiv\mathcal{U}_1$ of s-selfdecomposable was defined
by the non-linear shrinking operations (in short: s-operation) $U_r,
r>0$,\, (for $x>0,\, U_r(x):= \max(0, x-r)$); cf. Jurek (1981). [
Terminology: $"$s-selfdecomposable$"$ is abbreviation of that
$"$shrinking-selfdecomposable$"$.]

\begin{rem}
\emph{Since the process $Y$ has values in a metric separable
complete space  we may and do assume that the paths of $Y$ are
cadlag; cf. Theorem A.1.1 in Jurek and Vervaat (1983), p. 260. Since
the random integral in (1) is defined by a formal integration by
parts formula, therefore the random integral in question does exist;
cf. Jurek-Vervaat (1983), Lemma 1.1. \\
\indent Furthermore, since L\'evy processes are semi-martingales the
random integral (1) can be defined as the Ito stochastic integral.
However, for our purposes we do not need that generality. }
\end{rem}

For a positive natural $m$ and a sequence of positive real
$\beta_1,\beta_2,...,\beta_m$ and  a probability measure $\nu\in
ID$, let us define the mappings
\begin{equation*}
\mathcal{J}^{\{\beta_1,..., \beta_m\}}(\nu):=
\mathcal{J}^{\beta_m}(\mathcal{J}^{\{\beta_1,...,
\beta_{m-1}\}}(\nu))\,=\La\bl\int_{(0,1]}
t\;dY_{\mathcal{J}^{\{\beta_1,...,
\beta_{n-1}\}}(\nu)}(t^{\beta_m})\br\,.
\end{equation*}

Our main results say that the above composition can be written as a
single integral of the form $(\star \star)$ with a suitable chosen
time change $r$. Furthermore, the composition is expressed in terms
of the individual random integrals.

\begin{thm}
For positive reals $\beta_1,\beta_2,...,\beta_m$ and an infinitely
divisible probability measure  $\nu$ we have
\begin{equation}
\mathcal{J}^{\{\beta_1,..., \beta_m\}}(\nu)\,=\La\bl\int_{(0,1]}
t\;dY_{\nu}(r_{\{\beta_1,..., \beta_{m}\}}(t)\br=
I^{t,\,\,r_{\{\beta_1,..., \beta_{m}\}}}_{(0,1]}(\nu)
\end{equation}
and the time scale change $r_{\{\beta_1,..., \beta_{m}\}}$ is given
by
\begin{equation}
r_{\{\beta_1,..., \beta_{m}\}}(t)=
\mathbb{P}\,[\,U_1^{1/\beta_1}\cdot U_2^{1/\beta_2}\cdot ... \cdot
U_m^{1/\beta_m}\le t\,],  \ \ \ 0<t\le 1\,,
\end{equation}
where $U_i$'s are  mutually independent uniformly distributed, on
the unit interval, random variables.

If all $\be_1,...,\be_n$ are different then
\begin{equation}
r_{\{\beta_1,...,
\beta_{n}\}}(t):=\sum_{j=1}^n\,C_{j,n}\,t^{\beta_j}, \ \
C_{j,n}(\be_1,...,\be_n)\equiv C_{j,n}\,:=\prod_{k\neq j, k=1}^n\,
\frac{\beta_k}{\beta_k-\beta_j}
\end{equation}
and, in particular, we get the equality:  \quad $\sum_{j=1}^n
C_{j,n}=1$ .

If $\be_1= \be_2= ... = \be_m=\alpha\,\,(m\ge 1)$ then
\begin{equation}
r_{\{\underbrace{\alpha,...,\alpha}_{m-times}\}}(t)=
t^{\alpha}\,\sum_{j=0}^{m-1}\frac{(-\alpha\log t)^j}{j!} , \ \
\mbox{for} \ 0<t\le 1\,. \ \ \ \ \
\end{equation}
\end{thm}

\begin{rem}\emph{ (a) Note that if $\mathcal{E}(\lambda)$ denotes the exponential random
variable with the parameter $\lambda$ then
$e^{-\mathcal{E}(\lambda)} \stackrel{d}{=}U^{1/\lambda}$. Hence if
$\mathcal{E}_i(\alpha), 1\le i \le m,$ independent and identically
distributed exponential random variables then $r_{\{\alpha,...,
\alpha\}}$ is the cumulative distribution of $e^{-\mathcal{E}_1
(\alpha)}\cdot e^{-\mathcal{E}_2 (\alpha)}\cdot . . . \cdot
e^{-\mathcal{E}_m (\alpha)} \stackrel{d}{=} e^{- \gamma_{m, \alpha
}},$ where $\gamma_{m, \alpha}$ is the gamma random variable with
the shape parameter $m$ and the scale parameter $\alpha$, i.e., it
has the density of the form $\alpha^m /(m-1)!
\,\,x^{m-1}\,e^{-\alpha x}\,1_{(0,\infty)}(x)$; in particular, for
the case $\al=1$ see
Proposition 4 in Jurek (2004). \\
(b) For m=2, Theorem 1, and in particular the equality (4) were
shown in Czyzewska-Jankowska and Jurek (2009).}
\end{rem}

\medskip
\textbf{Example.} For $\beta_j:=j\,\beta$, $j=1,2,...,n$ with fixed
$\beta>0$, we get that
\[
C_{j,n}(\beta, 2\beta,...,n\beta)=(-1)^{j-1}{n \choose j} \ \mbox{
and} \ r_{\{\beta, 2\,\beta,...,
n\,\beta\}}(t):=\sum_{j=1}^n\,(-1)^{j-1}{n \choose j}\,t^{j\,\beta}.
\]

\medskip
In a proof of the above theorem the following identity, that might
be  also of an independent interest, is needed.
\begin{lem}
For distinct complex numbers $z_j, j=1,2,...,n,n+1$ we have
equality:
\begin{equation}
\sum^{n}_{i=1}\frac{1}{z_{i}-z_{n+1}}\, \Bigl(\prod^{n}_{k=1; k\neq
i}\frac{1}{z_{k}-z_{i}}\Bigr) =
\prod^{n}_{i=1}\frac{1}{z_{i}-z_{n+1}}.
\end{equation}
Equivalently, for any distinct complex numbers $z_j, j=1,2,...,n$,
we have the identity
\begin{equation}
\sum^{n}_{i=1}\,\prod^{n}_{k=1; k\neq i}\frac{z_k-z}{z_{k}-z_{i}}\,
\equiv\,1, \ \ \mbox{for all} \ \ \ z\in\Cset\,,
\end{equation}
that can be regarded as a decomposition of 1 as a sum of finite
products of complex fractions.
\end{lem}

\medskip
(It would be interesting to have a geometric description of the
above identity (7).)

\medskip
Since the characteristic function of each $\nu\in ID$  is uniquely
determined by the triple $[a, R, M]$ from its L\'evy-Khintchine
formula we will write formally that $\nu=[a,R,M]$; for details see
the Section 2.1 below.

If $\nu=[a,R,M]$ and $\mathcal{J}^{\{\be\}}(\nu):=[a^{\{\be\}},
R^{\{\be\}}, M^{\{\be\}}]$ and
\begin{equation}
b_{M,\be}:=\int_{\{||x||>1\}} x\,||x||^{-1-\be}\,M(dx)\in \Rset^d \
( \ \mbox{or a Banach space} \  E )
\end{equation}
then we have
\begin{multline}
a^{\{\be\}}= \be(\be+1)^{-1}(a+ b_{M, \be}), \ \ \
R^{\{\be\}}=\be(2+\be)^{-1}R, \\
M^{\{\be\}}(A)=\int_0^1\,T_{t^{1/\be}}M(A)dt = \be \int_0^1
M(s^{-1}A)\,s^{\be-1}ds, \ \mbox{for} \ A\in \mathcal{B}_0,
\end{multline}
where $\mathcal{B}_0$ stands for all Borel subsets of
$\Rset^d\setminus\{0\}$ (or $E\setminus\{0\}$. The above identities
follow from the  L\'evy-Khintchine formula (13) and Lemma 2; for
more details cf. Jurek (1988).

With these notations Theorem 1  gives the description of random
integrals (2) in terms of their corresponding triples.

\begin{thm}
For distinct positive reals $\beta_1,\beta_2,...,\beta_n$,
coefficients $C_{j,n}$ defined by (4), an infinitely divisible
probability measure $\nu=[a,R,M]$ and
$\mathcal{J}^{\{\beta_1,...,\beta_n\}}(\nu)=
[a^{\{\beta_1,...,\beta_n\}},R^{\{\beta_1,...,\beta_n\}},M^{\{\beta_1,...,\beta_n\}}]$
we have
\begin{gather}
a^{\{\beta_1,...,\beta_n\}}= a \prod_{j=1}^n \frac{\be_j}{\be_j+1}
+\sum_{j =1}^n \frac{\be_j\,b_{M,\be_j}}{\be_j+1}\prod_{k\neq j,
k=1}^n\,\frac{\beta_k}{\beta_k-\beta_j} =
 \sum_{j=1}^n C_{j,n}\,\,a^{\{\be_j\}} \\
R^{\{\beta_1,...,\beta_n\}}= \prod_{j=1}^n\frac{\be_j}{\be_j+2}\,R =
\sum_{j=1}^n C_{j,n} R^{\{\be_j\}}, \\
M^{\{\beta_1,...,\beta_n\}}(A)= \int_0^1...
\int_0^1\,T_{t_1^{1/\be_1}\,...\,t_n^{1/\be_n}}M(A)dt_1\,...\,dt_n
=\sum_{j=1}^n C_{j,n}\,M^{\{\be_j\}}(A)
\end{gather}
where $b_{M,\be_j}$, $a^{\{\be_j\}}$, $R^{\{\be_j\}}$ and
$M^{\{\be\}}$ are given in (8) and (9).
\end{thm}

\begin{cor}
For distinct positive reals $\beta_1,\beta_2,...,\beta_n$ and the
constants $C_{j,n}$ given in (4) we have
\[
\mathcal{J}^{\{\beta_1,...,\beta_n\}}(\nu)=(\mathcal{J}^{\{\be_1\}}(\nu))^{\ast
C_{1,n}} \ast ... \ast(\mathcal{J}^{\{\be_n\}}(\nu))^{\ast C_{n,n}},
\]
where for $C_{j,n}<0$ the corresponding convolution power means the
reciprocal of the corresponding infinitely divisible Fourier
transform.
\end{cor}

\textbf{2. Auxiliary results and proofs.}

\medskip
\textbf{2.1.} \underline{Random integrals.}

Let us recall that for a probability Borel measures $\mu$ on
$\Rset^d$ (or on $E$), its \emph{characteristic function} ( Fourier
transform) $\hat{\mu}$ is defined as
\[
\hat{\mu}(y):=\int_{\Rset^d} e^{i<y,x>}\mu(dx), \ y\in\Rset^d, \ \
(\mbox{or} \ \ y\in E^{\prime})
\]
where $<\cdot,\cdot>$ denotes the scalar product;  (in case of
Banach spaces, $<\cdot,\cdot>$ is the bilinear form on
$E^{\prime}\times E$). Further, the characteristic  function of an
infinitely divisible probability measure $\mu$  admits the following
\emph{L\'evy-Khintchine representation:}
\begin{multline}
\hat{\mu}(y)= e^{\Phi(y)}, \ y \in \Rset^d, \ \ \mbox{and the L\'evy
exponent} \ \ \Phi(y)=i<y,a>- \\ \frac{1}{2}<y,Ry> + \int_{\Rset^d
\backslash \{0 \}}[e^{i<y,x>}-1-i<y,x>1_B(x)]M(dx),  \qquad
\end{multline}
where $a$ is  a \emph{shift vector}, $R$ is a \emph{covariance
operator} corresponding to the Gaussian part of $\mu$, $B:=\{x:
||x||\le1\}$ (the unit ball) and $M$ is a \emph{L\'evy spectral
measure}. (We add the term \emph{spectral} to avoid the possible
confusion with \emph{L\'evy measures} as sometimes are called the
selfdecomposable (class L) measures.) Since there is a one-to-one
correspondence between measures $\mu \in ID$ and triples $a$, $R$
and $M$ in its L\'evy-Khintchine formula (13) we will formally write
$\mu=[a,R,M]$.

Note that for $s\in\Rset$ we have
\begin{multline}
\Phi(sy)=i<y, s(a+\int_{E\setminus\{0\}}x(1_B(sx)-1_B(x))M(dx))>
-\frac{1}{2}s^2<y,Ry> \\ + \int_{E \backslash \{0
\}}[e^{i<y,z>}-1-i<y,z>1_B(z)]M(s^{-1}dz)
\end{multline}
Finally, let us recall that
\begin{equation}
M \ \mbox{is L\'evy spectral measure on $\Rset^d$ iff} \ \
\int_{\Rset^d}\min(1, ||x||^2)M(dx)<\infty
\end{equation}
For infinity divisibility on Banach spaces we refer to the monograph
by Araujo and Gin\'e (1980) or Linde (1983) or Ledoux and Talagrand
(1991). Let us emphasize here that the characterization (15) of
L\'evy spectral measures is NOT true on infinite dimensional Banach
spaces ! However, it holds true on Hilbert spaces; cf. Parthasarathy
(1967), Chapter VI.

\medskip
For this note it is important to have the following technical
result:
\begin{lem}
If the random integral $A\equiv \int_{(a,b]} h(t)dY_{\nu}(r(t))$
exists then we have
\begin{equation*}
\log \widehat{\mathcal{L}(A)}(y)=\int_{(a,b]} \log
\widehat{\mathcal{L}(Y_{\nu}(1))}(h(s)y)dr(s)=\int_{(a,b]}\,\Phi(h(s)y)dr(s),
\end{equation*}
where $y\in \Rset^d$ (or $E^{\prime}$) and $\Phi$ is the L\'evy
exponent of $\widehat{\mathcal{L}(Y_{\nu}(1))}=\hat{\nu}$. In
particular, if $r$ is the cumulative probability distribution
function of a random variable $T$ concentrated of the interval
$(a,b]$ then  \ $\log
\widehat{\mathcal{L}(A)}(y)=\mathbb{E}[\Phi(h(T)y)]$.
\end{lem}
The formula in Lemma 2 is a straightforward consequence of our
definition (integration by parts)  of the random integrals ($\star
\star$). The proof is analogous to that in  Jurek-Vervaat(1983),
Lemma 1.1 or Jurek (1988), Lemma 3.2 (b).

\begin{rem}
\emph{Note that for  bounded intervals $(a,b] \subset\Rset^+$,
positive monotone functions $r$  and real-valued continuous, bounded
variation functions $h$, the integrals of the form $A$ in Lemma 2
 are well-defined.}
\end{rem}

\medskip
\medskip
\textbf{2.2.} \underline{Proof of Lemma 1.} For each $1\le i\le n$,
the polynomials $Q_i$, of $n-1$ degree,  given by the following
\begin{equation}
Q_i(z):=\prod^{n}_{k=1; k\neq i}\frac{z_k-z}{z_{k}-z_{i}} \ \
\end{equation}
satisfy the conditions
\[ Q_i(z_i)=1, \ \ Q_i(z_j)=0 \ \mbox{for} \ \ 1\le i\neq j \le n\,.
\]
Consequently, the function
\begin{equation}
\textbf{Q}(z):=\sum_{i=1}^n\,Q_i(z)-1 \ \ \mbox{is vanishing in n
points} \ z_1,z_2,...,z_n.
\end{equation}
Since \textbf{Q} is a polynomial of n-1 degree we conclude that
$\textbf{Q}(z)\equiv 0$ which completes the proof of Lemma 1.

\medskip
\begin{rem}
It might be of an additional interest to recall here that in the
interpolation theory for given set of points $(x_0,y_0),$
$(x_1,y_1), ..., (x_n,y_n)$ in the plane $\Rset^2$, with distinct
$x_0,x_1,...,x_n$,
\[
P_n(x):=\sum_{j=0}^n\,y_j\prod_{k\neq j,
k=0}^n\frac{x_k-x}{x_k-x_j},
\]
is the unique Lagrange interpolating polynomial of degree less or
equal $n-1$ and such that
\[
P_n(x_i)=y_i \ \ \mbox{for all} \ \ i=0,1,2,...,n;
\]
cf.  Kincaid and Cheney (1996), Chapter 6. Thus for the particular
points $(x_0,1), (x_1,1),..., (x_n,1)$ in $\Rset^2$ we get  the line
$y=P_n(x)=1$ as the Lagrange interpolating polynomial.
\end{rem}

\medskip
\medskip
\textbf{2.3.} \underline{Products of independent uniformly
distributed random variables.}

Here are some elementary identities concerning the products of
powers of independent uniformly distributed  random variables. The
main objective is to express the cumulative distribution function
(c.d.f.) or probability density function (p.d.f.) of such products
as a linear combinations (with not necessary positive coefficients)
of other c.d.f. (or p.d.f.).
\begin{lem}
Let $U_i$, $1\le i \le n$ be i.i.d. uniformly distributed over the
interval $(0,1]$,\ $\alpha_i>0$  and let $f_{\{\alpha_1,
\alpha_2,...,\alpha_n\}}$ and $F_{\{\alpha_1,
\alpha_2,...,\alpha_n\}}$ denote the probability density and
cumulative distribution function of $U^{1/\alpha_1}_1\cdot
U_2^{1/\alpha_2}\cdot ...\cdot U^{1/\alpha_n}_n$, respectively. Then
\begin{multline}
(a) \ \ f_{\{\alpha_1, \alpha_2,...,\alpha_n\}}(x_n)=
\\ \alpha_n\,x_n^{\alpha_n-1}\int_{x_n}^1
\alpha_{n-1}\,x_{n-1}^{\alpha_{n-1}-\alpha_n-1} \int_{x_{n-1}}^1
\alpha_{n-2}\,x_{n-2}^{\alpha_{n-2}-\alpha_{n-1}-1}
\int_{x_{n-2}}^1\alpha_{n-3}\,x_{n-3}^{\alpha_{n-3}-\alpha_{n-2}-1} ... \\
...\int_{x_3}^1 \alpha_2 \,x_2^{\alpha_2-\alpha_3-1}\int_{x_2}^1
\alpha_1 \,x_1^{\alpha_1-\alpha_2-1}\,dx_1\,dx_2
\,...\,dx_{n-2}\,dx_{n-1}, \ 0<x_n \le 1.
\end{multline}

\noindent (b) If $\alpha_1= \alpha_2= ... = \alpha_m=\alpha$ then
\begin{gather}
f_{\{\underbrace{\alpha,...,\alpha}_{m-times}\}}(x)= \alpha\,
x^{\alpha
-1}\,\frac{(-\alpha\,\log x )^{m-1}}{(m-1)!} \ \ \ \mbox{for} \ \ 0<x\le 1, \\
F_{\{\underbrace{\alpha,...,\alpha}_{m-times}\}}(s)=
s^{\alpha}\,\sum_{j=0}^{m-1}\frac{(-\alpha\log s)^j}{j!} \  \ \
\mbox{for} \ \ 0<s\le 1, \ \ \ \ \
\end{gather}
and $F_{\{\underbrace{\alpha,...,\alpha}_{m-times}\}}(s)= 1$ for
$s\ge1$ and zero for $s<0$.

\noindent (c) If all positive reals $\alpha_i$ ($i=1,2,...,n$) are
distinct and
\begin{equation}
C_{j,n}(\al_1,...,\al_n)\equiv C_{j,n}:=\prod_{k \neq j\, ,
k=1}^n\frac{\alpha_k}{\alpha_k-\alpha_j} \ \ \ \mbox{and} \ \
c_{j,n}:=\prod_{k \neq j\, , k=1}^n\frac{1}{\alpha_k-\alpha_j}
\end{equation}
then
\begin{equation}
f_{\{\alpha_1, \alpha_2,...,\alpha_n\}}(x)=
\sum_{j=1}^n\,C_{j,\,n}\,\alpha_j
x^{\al_j-1}=\al_1...\al_n\sum_{j=1}^n c_{j,n}x^{\al_j-1}, \  \
0<x\le 1,
\end{equation}
and
\begin{equation}
F_{\{\alpha_1,
\alpha_2,...,\alpha_n\}}(s)=\sum_{j=1}^n\,C_{j,\,n}\,\,s^{\alpha_j}
,\ \mbox{for} \ 0<s\le 1, \ \ \mbox{where} \
\sum_{j=1}^n\,C_{j,\,n}=1.
\end{equation}
\end{lem}

\emph{Proof.} For positive and independent rv $X$ and $Z$ with
p.d.f. $f_X$ and $f_Z$, respectively we have  that $X\cdot Z$ has
the p.d.f.
\begin{equation}
f_{X\cdot Z}(z)=\int_0^{\infty}f_X(\frac{z}{x})\frac{1}{x}f_Z(x)\,
dx \, .
\end{equation}
Since $f_{\{\alpha\}}(x)=\alpha\,x^{\alpha-1}\,1_{(0,1)}(x)$ is the
p.d.f. of $U^{1/\alpha}$ therefore from (24) we get
\begin{equation}
f_{\{\al , \be \}}(z)= f_{U^{1/\al} \cdot U^{1/\be}}(z) =
\be\int_z^1f_{\{\al\}}(\frac{z}{x}) x^{\be-2}dx
 =f_{\{\be\}}(z)\int_z^1 f_{\{\al\}}(t)\,t^{-\be}\,dt.
\end{equation}
Hence for  $\alpha_1$ and $\alpha_2$  we conclude that
\begin{gather}
f_{\{\alpha_1,\alpha_2\}}(x_2)=\alpha_2\,x_2^{\alpha_2-1}\int_{x_2}^1\alpha_1\,x_1^{\alpha_1-\alpha_2-1}\,dx_1
\end{gather}
which is indeed of the form (18) for $n=2$.

Assume,  by the mathematical induction argument,  that the formula
(18) holds true for $n$. Then using (24) and (25) we obtain
\begin{multline*}
f_{\{\alpha_1, \alpha_2,...,\alpha_n,
\alpha_{n+1}\}}(x_{n+1})=f_{\alpha_{n+1}}(x_{n+1})\int_{x_{n+1}}^1
f_{\{\alpha_1, \alpha_2,...,\alpha_n\}}(x_n)\,
x_n^{-\alpha_{n+1}}dx_{n} = \\
\alpha_{n+1}x_{n+1}^{\alpha_{n+1}-1}\int_{x_{n+1}}^1\,\Big(\alpha_n\,x_n^{\alpha_n-1}\int_{x_n}^1
\alpha_{n-1}\,x_{n-1}^{\alpha_{n-1}-\alpha_n-1} \int_{x_{n-1}}^1
\alpha_{n-2}\,x_{n-2}^{\alpha_{n-2}-\alpha_{n-1}-1} \\ ...
\int_{x_3}^1 \alpha_2 \,x_2^{\alpha_2-\alpha_3-1}\int_{x_2}^1
\alpha_1 \,x_1^{\alpha_1-\alpha_2-1}\,dx_1\,dx_2
\,...\,dx_{n-2}\,dx_{n-1}\Big) \ x_n^{-\alpha_{n+1}}dx_{n}= \\
\alpha_{n+1}x_{n+1}^{\alpha_{n+1}-1}\int_{x_{n+1}}^1\,\alpha_n\,x_n^{\alpha_n-\alpha_{n+1}-1}\int_{x_n}^1
\alpha_{n-1}\,x_{n-1}^{\alpha_{n-1}-\alpha_n-1} \int_{x_{n-1}}^1
\alpha_{n-2}\,x_{n-2}^{\alpha_{n-2}-\alpha_{n-1}-1}\\ ...
\int_{x_3}^1 \alpha_2 \,x_2^{\alpha_2-\alpha_3-1}\int_{x_2}^1
\alpha_1 \,x_1^{\alpha_1-\alpha_2-1}\,dx_1\,dx_2
\,...\,dx_{n-1}\,dx_{n},
\end{multline*}
which is  the equality (18) for $n+1$. Thus the proof of the part
(a) is complete.

Taking in (18),  $\al_1=\al_2=...=\al_n=\al$ and performing the
successive integrations, we get the formula (19). (Or simply prove
(19) by the induction argument utilizing (24)). Integrating p.d.f.
(19) we get c.d.f. (20) and this establishes the part (b).

In the part(c), formulae (22) and (23) are obvious for n=1. Assume
that (22) holds true for n.  First, from Lemma 1 formula (6) we
infer that for $1\le j \le n$ we get
\begin{equation}
c_{j,n}(\alpha_{n+1}-\alpha_{j})^{-1}=  c_{j,n+1},\ \
 \ \  \sum_{j=1}^n\,\
(\alpha_j-\alpha_{n+1})^{-1}c_{j,n}=  c_{n+1, n+1} \ .
\end{equation}
Then from (25), (22) and (27) and again (6) from Lemma 1 we get
\begin{multline*}
f_{\{\alpha_1, \alpha_2,...,\alpha_n,
\alpha_{n+1}\}}(x)=\alpha_{n+1}x^{\alpha_{n+1}-1}\int_{x}^1
f_{\{\alpha_1, \alpha_2,...,\alpha_n\}}(t)\, t^{-\alpha_{n+1}}dt \\
=\alpha_1...\alpha_{n+1}\sum_{j=1}^n\,c_{j,n}\,\,x^{\alpha_{n+1}-1}\int_{x}^1
t^{\alpha_j-\alpha_{n+1}-1}dt\\
=\alpha_1...\alpha_{n+1}\sum_{j=1}^n\,c_{j,n}\frac{1}{\alpha_j-\alpha_{n+1}}
(x^{\alpha_{n+1}-1}-x^{\alpha_j-1})\\
=\alpha_1...\alpha_{n+1}\sum_{j=1}^n\,c_{j,n+1}x^{\alpha_j-1}+
\alpha_1...\alpha_{n+1}\big(\sum_{j=1}^n\,c_{j,n}\frac{1}{\alpha_{j}-\alpha_{n+1}}\big)x^{\alpha_{n+1}-1}\\
= \alpha_1...\alpha_{n+1} (  \sum_{j=1}^n\,c_{j,n+1}x^{\alpha_j-1} +
c_{n+1,n+1}\,\, x^{\al_{n+1}-1} ) =
\alpha_1...\alpha_{n+1}\,\sum_{j=1}^{n+1}\,c_{j,n+1}x^{\alpha_j-1},
\end{multline*}
which completes the proof of (22) and consequently of (23). This
completes the proof of Lemma 3.

\medskip
\medskip
\textbf{2.4.} \underline{Proof of Theorem 1.}

In view of Lemma 2, to prove formula (2) it is necessary and
sufficient to show the equality
\begin{multline}
\log\bl(\mathcal{J}^{\{\beta_1,...,
\beta_n\}}(\nu))^{\wh}\br(y)=\mathbb{E}[\log\hat{\nu}(U_1^{1/\be_1}\cdot
U_2^{1/\be_2}\cdot...\cdot U_n^{1/\be_n}y)]
\\ =\int_0^1\log\hat{\nu}(ty)dr_{\{\beta_1,..., \beta_{n}\}}(t) ,  \ \ y \in
E^{'}. \ \ \ \ \ \ \
\end{multline}
Of course, (28) holds for $n=1$. Assume it is true for $n-1$. Then
from Lemma 2 and the definition of the  mapping
$\mathcal{J}^{\{\be_1,...,\be_n\}}$ (given before Theorem 1) we get
\begin{multline}
\log(\mathcal{J}^{\{\beta_1,..., \be_{n-1} \beta_n\}}(\nu))^{\wh}(y)
= \int_0^1 \log(\mathcal{J}^{\{\beta_1,...,
\beta_{n-1}\}}(\nu))^{\wh}(s\,y)ds^{\be_n} = \\ \int_0^1
\mathbb{E}[\log\hat{\nu}(U_1^{1/\be_1}\cdot
U_2^{1/\be_2}\cdot...\cdot U_{n-1}^{1/\be_{n-1}}s y)]ds^{\be_n} =
\mathbb{E}[\log\hat{\nu}(U_1^{1/\be_1}\cdot
U_2^{1/\be_2}\cdot...\cdot U_n^{1/\be_n}y)]\\ =
\int_0^1\log\hat{\nu}(sy)dr_{\{\beta_1,...,
\beta_{n-1},\beta_n\}}(t), \ \ \ \ \ \ \ \
\end{multline}
which completes  proof of (28) and consequently  the formulae (2)
and (3). Explicit expressions for time changes $r_{\{\beta_1,...,
\beta_{n-1},\beta_n\}}(t)$ are given in Lemma 3, part (c).

\medskip
\medskip
\medskip
\textbf{2.5.} \underline{Proof of Theorem 2.}

First, putting $\Phi(y)= \log\hat{\nu}(y)$ ( the L\'evy exponents of
$\nu$) into (28) we get
\begin{multline*}
i<y,a ^{\{\be_1,...,\be_n\}}>-  \frac{1}{2}<y,R^{\{\be_1,...,\be_n\}}y> \\
 + \int_{\Rset^d
\backslash \{0
\}}[e^{i<y,x>}-1-i<y,x>1_B(x)]M^{\{\be_1,...,\be_n\}}(dx)\\
=\mathbb{E}[\Phi(U_1^{1/\be_1}\cdot U_2^{1/\be_2}\cdot...\cdot
U_n^{1/\be_n}y)] \ \
\end{multline*}
Since for $0<s\le 1$, we have that $1_B(sx)-1_B(x)=1_{\{1<||x||\le
s^{-1}\}}(x)$ therefore from the above and (14) we get
\begin{multline*}
a^{\{\beta_1,...,\beta_n\}}=\mathbb{E}[U_1^{1/\be_1}\cdot
U_2^{1/\be_2}\cdot...\cdot U_n^{1/\be_n}(a+
\int_{1<||x||\le(U_1^{1/\be_1}\cdot U_2^{1/\be_2}\cdot...\cdot
U_n^{1/\be_n})^{-1}}\,x M(dx))]\\
=\prod_{j=1}^n\frac{\be_j}{1+\be_j}\, a +\int_0^1\int_{1<||x||\le
s^{-1}}\,s x\,M(dx)dr_{\{\beta_1,..., \beta_{n}\}}(s) \ \ \ \
(\mbox{by} \ (3) \ \mbox{and}\ (4)) \\
=\prod_{j=1}^n\frac{\be_j}{1+\be_j}\,a  + \be_1\be_2...\be_n
\sum_{j=1}^n\,c_{j,n}\int_0^1\int_{1<||x||\le s^{-1}}\,s
x\,M(dx)s^{\be_j-1}ds   \ \  (\mbox{by (8)}) \\
=\prod_{j=1}^n\frac{\be_j}{1+\be_j}\, a  +\sum_{j =1}^n
\frac{\be_j}{\be_j+1}\bl\prod_{k\neq j,
k=1}^n\,\frac{\beta_k}{\beta_k-\beta_j}\br\,b_{M,\be_j} \ \
(\mbox{by (7) with z = -1})\\=
\prod_{j=1}^n\frac{\be_j}{1+\be_j}\,[\sum_{j=1}^n(a+b_{M,\be_j})\prod_{k\neq
j, k=1}^n\frac{\be_k+1}{\be_k-\be_j}]=\sum_{j=1}^n
\,C_{j,n}\,a^{\{\be_j\}},
\end{multline*}
which proves the formula for the shift vector.

Similarly, for the Gaussian part, using again the identity  (7)
(with $z_j=\be_j$ and $z=-2$) we get
\begin{multline*}
R^{\{\beta_1,...,\beta_n\}}=\mathbb{E}[(U_1^{1/\be_1}\cdot
U_2^{1/\be_2}\cdot...\cdot
U_n^{1/\be_n})^2]\,R=\prod_{j=1}^n\frac{\be_j}{\be_j+2}\,R \\
=\prod_{j=1}^n\frac{\be_j}{\be_j+2}\bl \sum_{l=1}^n\prod_{k\neq l,
k=1}^n\frac{\be_k+2}{\be_k-\be_l}\br\,R=\sum_{l=1}^n(\prod_{k\neq l,
k=1}^n\frac{\be_k}{\be_k-\be_l})\,R^{\{\be_l\}}= \sum_{l=1}^n
C_{l,n}\,R^{\{\be_l\}},
\end{multline*}
which gives the formula for Gaussian covariance.

Finally for the L\'evy spectral measure using (28), (4) and (9) we
have
\begin{multline*}
M^{\{\beta_1,...,\beta_n\}}(A)=\mathbb{E}[M((U_1^{1/\be_1}\cdot
U_2^{1/\be_2}\cdot...\cdot
U_n^{1/\be_n})^{-1}A)]\\
=\int_0^1\,T_{s}M(A)dr_{\{\beta_1,...,
\beta_{n}\}}(s)=\be_1\be_2...\be_n\,
\sum_{j=1}^n\,c_{j,\,n}\,\int_0^1M(s^{-1}A)\,s^{\be_j-1}ds \\
 =\sum_{j=1}^n (\prod_{k\neq j, k=1}^n
 \frac{\be_k}{\be_k-\be_j})\int_0^1\,M(s^{-1/\be_j}A)\,ds =
\sum_{j=1}^n \, C_{j,n}\,M^{\{\be_j\}}(A),
\end{multline*}
which completes the proof of Theorem 2.

\medskip
\textbf{3. Concluding remarks.} Although generalized
s-selfdecomposable distributions are known for $\beta<0$ analogous
calculus on the corresponding random integrals seems to be more
complicated. Similarly, Theorem 1  for not necessarily all different
$\beta's$ might be much more involved and is not discussed in this
note.

\medskip
\textbf{Acknowledgements.} Author would like to thank the Reviewer
for suggesting another proof for Lemma 1. Original proof used the
mathematical induction arguments.

\medskip
\begin{center}
\textbf{References}
\end{center}

\noindent [1] D. J. Aldous (1977), Limit theorems for subsequences
of arbitrary-dependent sequences of random variables, \emph{Z.
Wahrscheinlichkeitstheorie verw. Gebiete},  \textbf{40}, pp. 59 --
82.

\noindent [2] A. Araujo  and E. Gine (1980), \emph{The central limit
theorem for real and Banach valued random variables.} John Wiley \&
Sons, New York.

\noindent[3] A. Czyzewska-Jankowska and Z. J. Jurek (2009), A note
on a composition of two random integral mappings $\mathcal{J}^{\be}$
and some examples, \emph{Stochastic Anal. Appl. } vol. 27 No 6, pp.
1212 -1222.

\noindent[4] L. James and Z. Zhang (2012), Quantile clocks,
submitted to \emph{Ann. Appl. Probab.}

\noindent [5] Z. J. Jurek (1981), Limit distributions for sums of
shrunken random variables, \emph{Dissertationes Mathematicae}, vol.
CLXXXV, PWN Warszawa .

\noindent [6] Z. J. Jurek (1985), Relations between the
s-selfdecomposable and selfdecomposable measures. \emph{Ann.
Probab.} vol. 13, Nr 2, pp. 592-608.

\noindent [7] Z. J. Jurek (1988), Random integral representation for
classes of limit distributions similar to L\'evy class $L_{0}$,
 \emph{Probab. Th. Rel. Fields.}  vol. 78, pp. 473-490.

\noindent [8] Z. J. Jurek (1989), Random Integral representation for
Classes of Limit Distributions Similar to L\'evy Class $L_{0}$, II ,
\emph{Nagoya Math. J.}, vol. 114, pp. 53-64.

\noindent [9] Z. J. Jurek (2004), The random integral representation
hypothesis revisited: new classes of s-selfdecomposable laws. In:
Abstract and Applied Analysis; \emph{Proc. International Conf.
ICAAA, Hanoi, August 2002}, World Scientific, Hongkong, pp. 495-514.

\noindent[10] Z. J. Jurek (2011). The random integral representation
conjecture: a quarter of a century later. \emph{Lithuanian Math.
Journal} vol. 51, no 3, pp. 362-369.

\noindent[11]  Z. J. Jurek and B. M. Schreiber (1992), Fourier
transforms of measures from classes $\mathcal{U}_{\be}$, $-2<\be\le
1$, \emph{J. Multivar. Analysis}, vol. 41, pp.194-211.

\noindent[12]  Z. J. Jurek and W. Vervaat (1983), An integral
representation for selfdecomposable Banach space valued random
variables, \textit{Z. Wahrsch. verw. Gebiete}, vol. 62, pp. 247-262.

\noindent [13] D. Kincaid and W. Cheney (1996), \emph{Numerical
Analysis; Mathematics of Scientific Computing}, Brooks/Cole
Publishing Company, Second Edition.

\noindent[14] M. Ledoux and M. Talagrand (1991).  \emph{Probability
in Banach spaces}, Springer-Verlag.

\noindent[15]  W. Linde (1983), \emph{Probability measures in Banach
spaces- stable and infinitely divisible distributions}, Wiley, New
York.

\noindent[16] K. R. Parthasarathy (1967), \emph{Probability measures
on metric spaces}. Academic Press, New York and London.

\medskip
\medskip
\medskip
\noindent Institute of Mathematics, University of Wroc\l aw \\ Pl.
Grunwaldzki 2/4, 50-384 Wroc\l aw, Poland \\
e-mail: zjjurek@math.uni.wroc.pl \ \ \
www.math.uni.wroc.pl/$\sim$zjjurek

\end{document}